\documentclass[12pt]{amsart}
\usepackage{graphicx}
\usepackage{xcolor}
\usepackage{hyperref}
\usepackage{fullpage}
\usepackage{listings}

\begin{document}

\title{Using the HOMFLY-PT polynomial to compute knot types}

\author{Eric J. Rawdon}
\address{Department of Mathematics, University of St.~Thomas, Saint Paul, MN 55105, USA}
\email{ejrawdon@stthomas.edu}
\urladdr{george.math.stthomas.edu/$\sim$rawdon}

\author{Robert G. Scharein}
\address{Hypnagogic Software, Vancouver, BC, V6K 1V6, Canada}
\email{rob@hypnagogic.net}
\urladdr{www.knotplot.com}




 
\begin{abstract}
  The HOMFLY-PT polynomial is a link invariant which is effective in
  determining chiral knot and link types with small crossing numbers.
  In this chapter, we concentrate on knots.  We provide a
  guide for computing the knot types of configurations from 3D
  coordinates via the HOMFLY-PT polynomial using publicly-available
  Linux freeware.  We include data on the efficacy of HOMFLY-PT for
  knot types through crossing number 16.
\end{abstract}
 
\maketitle

\section{Introduction}

The goal of this article is to provide a guide for researchers to
compute knot types from vertex coordinates of polygons or knot
crossing information using the HOMFLY-PT polynomial.  We will
concentrate on software that is freely available and note potential
issues.  The software can also be used to compute link types, although
the HOMFLY-PT polynomial has limitations which make it less effective
when working with certain classes of links, the details of which are
explained below.

The tools described here are designed for large-scale computations,
e.g., random knot studies generating thousands or millions of knots,
although they can also be used for a single knot configuration.


The HOMFLY-PT polynomial is a two-variable Laurent polynomial on
oriented links which was discovered by four different research
groups \cite{HOMFLY,PrzTra}. It is a generalization of the
Alexander \cite{Alexander} and Jones \cite{Jones} polynomials.  There
are different normalizations for the variables.  We use the
normalization proposed by Lickorish and Millett \cite{LicMil}, which
we call the \textit{LM-normalization}.

The HOMFLY-PT polynomial cannot distinguish all knot types.  There are
distinct knot types sharing the same HOMFLY-PT polynomial, which
we will call \emph{collisions}, for even small crossing number knot
types (e.g., $+5_1$ and $+10_{132}$ share the same HOMFLY-PT
polynomial).  Furthermore, knot mutants have the same HOMFLY-PT
polynomial \cite{Adams}.  
It is unknown if the only knot type with a HOMFLY-PT polynomial of 1
is the unknot (although none have yet been found).  These facts might
make the user pause.  However, in our opinion, the HOMFLY-PT
polynomial provides a good balance of identifying knot types in a
reasonable amount of time for knot types through crossing number 16.

Since the HOMFLY-PT polynomial is a generalization of the Jones
polynomial \cite{Jones}, which, in turn, is a generalization of the
Alexander polynomial, the number of collisions for HOMFLY-PT is
smaller than the number of collisions for these other polynomials.
Table \ref{collisiontable2} shows the number of collisions accumulated
through different crossing numbers.  One can find all of the
collisions in the file \texttt{jhomflytable.txt} (available at
\cite{knottingtools} and described in Section \ref{sec:prefixes}) for
knot types through crossing number 16.  Other techniques can be used
when there are collisions, which we briefly explain in Section
\ref{sec:collisions}.

The HOMFLY-PT polynomial computation is exponential time on the number of
crossings \cite{Vertigan}.  
The Alexander polynomial can be computed in $O(c^3)$ time, where $c$
is the number of crossings in the given knot diagram, according to the
algorithm given in Alexander's original paper \cite{Alexander} (it
involves computing a matrix determinant of small polynomials).  If one
is analyzing a set of knots where the configurations have many
crossings, but the knot types of these knots have small crossing
number, the Alexander polynomial could be a better choice.

The HOMFLY-PT polynomial is an invariant of oriented links, although
it is not sensitive to orientation on knots (and really only relative
orientation on links).  For the sake of knots, the HOMFLY-PT
polynomial works at the level of pictures of knot types.
Specifically, given an oriented knot configuration $K$ of a knot type
$\mathcal{K}$, let $K_r$ denote $K$ with the opposite orientation,
$K^*$ denote the mirror image of $K$, and $H(K)$ be the HOMFLY-PT
polynomial of $K$.  Then for every knot configuration $K$, it is the case that
$$H(K) = H(K_r)\,.$$
If $K$ is isotopic to either $K^*$ or $K^*_r$ (or both), then
the knot type $\mathcal{K}$
is said to be \textit{achiral}\footnote{Note that achiral knot types are
  also called \textit{amphichiral} or \textit{amphicheiral}.}.  If $K$ is isotopic to neither $K^*$
nor $K^*_r$, then the knot type is said to be \textit{chiral}.  If $K$ is
from a chiral knot type, it is usually, but not always, the case that
$H(K)\not = H(K^*)$.

The simplest knot type where $K$ and $K_r$ are not isotopic is
$8_{17}$ \cite{Adams}.  Interestingly, $8_{17}$ is equivalent to its
mirror image with the opposite orientation.  We invite
readers to search for terms like \textit{invertible knot},
\textit{fully achiral knot}, \textit{positive achiral knot}, and
\textit{negative achiral knot} to further understand the symmetry
groups of knot types, or consult Chapter 10 of Kawauchi's
book~\cite{Kawauchi}.

Unfortunately, there do not appear to be simple tools to distinguish
between oriented knot types, which is unfortunate since it could be of
interest in scientific applications.  Note that HOMFLY-PT is sensitive
to orientations for links, but only relative orientation.  For
example, there are four possible choices for orientations for a
2-component link $L$.  However, there are at most two different
HOMFLY-PT polynomials amongst these four orientation choices for $L$.

The HOMFLY-PT polynomial can (often) be used to detect the chirality of a knot type.
For example, there are two chiralities for the trefoil knot, the
so-called right- and left-handed trefoils (which we denote $+3_1$ and
$-3_1$, respectively), whose HOMFLY-PT polynomials differ.  However,
there are many knot types which are chiral, yet both chiralities share
the same HOMFLY-PT polynomial ($9_{42}$ is the knot type with lowest
crossing number with this property).

Also note that mathematicians typically refer to knot types without
specifying a chirality, for example saying $3_1$ instead of specifying
$+3_1$ versus $-3_1$.  We refer to, for example, $3_1$ as a
\textit{base knot type} and $+3_1$ as a \textit{chiral knot type}.
For scientific applications, the chirality can be critical, which is
why we specify chirality when possible.  See Section
\ref{sec:chiralknottypes} for data about chiral knot types.

Note that until recently, prime base knot types had been fully classified
up through crossing number 16 \cite{HosThiWee}.  In 2020, Burton
\cite{Burton} pushed the enumeration through crossing number 19, resulting
in 352,152,252 non-trivial prime base knot types.  We
are not aware of tables of chiral knot types (although the data
clearly has been available for a long while) so have included
counts for the number of chiral knot types by crossing number
through crossing number 16 in Section \ref{sec:chiralknottypes}.

The numbers above are only for prime knot types.  One can create
composite knot types by combining two or more non-trivial prime knot
types \cite{Adams}.  It has not been shown that the crossing number is
additive over connected sums (i.e., knot composition).  However, for
the sake of this paper, we will act as though the crossing number of a
$p$-crossing prime knot type and a $q$-crossing prime knot type is
$p+q$, in which case, we can fully enumerate all composite knot types
through 19 crossings.  However, we consider only prime and composite
knot types whose crossing number does not exceed 16, amongst which
there are 3,473,712 chiral knot types and 2,427,782 unique HOMFLY-PT
polynomials.

Most of the computations discussed here can also be applied to links.
Prime links have been enumerated through 16 crossings
\cite{hostelinkenumeration} for up to 4-component links.  Link
composition is more complicated than knot composition because one
needs to specify onto which component the composition is completed.
Furthermore, the HOMFLY-PT polynomial is multiplicative over
compositions, so $H(L_1\#L_2)=H(L_1)\cdot H(L_2)$.  Due to this
multiplicative nature of HOMFLY-PT, the classification of composite
links using the HOMFLY-PT polynomial is less effective than for knots.

There are several packages which can compute the HOMFLY-PT polynomial
from some set of input data.  These include the Ewing-Millett code
\texttt{lmpoly} \cite{millettewing}, code (which is not public) by
Gouesbet et al.~\cite{ghomfly}, \texttt{topoly} by Su{\l}kowska's
group \cite{topoly}, the command homfly\_polynomial in \texttt{Sage} \cite{sage}, the command HOMFLYPT (written by Scott Morrison)
in the KnotTheory package \cite{knottheorypackage} for
\texttt{Mathematica}, and modifications of Jenkins's code
\cite{Jenkins,libhomfly,knottingtools}, which is the focus of this
paper.

\section{Using freeware to compute HOMFLY-PT polynomials and knot types}

The software described in this section is freeware and can be
downloaded from \cite{knottingtools}.  The user must also install
\texttt{libhomfly}, which is available for easy installation from
several Linux repositories as well as from GitHub \cite{libhomfly}.
We describe the installation process briefly in Section
\ref{sec:install} and more detailed information can be found at
\cite{knottingtools}.

\subsection{Overview of the process}
\label{sec:overview}
We begin at the stage of having a list of 3D coordinates forming a
polygonal knot.  The computation proceeds as
$$\text{coordinates} \rightarrow \text{EGCs} {\color{gray}\rightarrow \text{simplified EGCs}}\rightarrow
\text{HOMFLY-PT} \rightarrow \textbf{knot type}$$ where 
\textbf{knot type} is the set of knot types with crossing number 16
or smaller matching the given HOMFLY-PT polynomial.  An EGC
is an extended Gauss code, which is one way to encode
the crossing data from a knot or link projection.  More information
about EGCs can be found in Section \ref{sec:egc}.  The grayed step to
simplify the EGCs is optional.  Also note that a user can begin at any
of the stages if they have appropriately formatted files.

For each stage in the above procedure, we provide open source public
domain software, which can be downloaded from \cite{knottingtools}.  In
addition, the software KnotPlot \cite{KnotPlot} can be used to compute
crossing codes, which has the same code base as \texttt{coords2egc}. In
particular, below is a list of these steps and the software used for
these tasks.

\begin{itemize}
\item \texttt{coords2egc}: $\text{coordinates} \rightarrow \text{EGC code}$
\item \texttt{xinger}: $\text{EGCs} \rightarrow{\color{gray} \text{simplified EGCs}}$
\item \texttt{jhomfly}: $\text{EGCs or {\color{gray}simplified EGCs}}\rightarrow \text{HOMFLY-PT}$
\item \texttt{jidknot}, \texttt{jidknot\_table16.py}, or \texttt{jidknot\_grep.py}: $\text{HOMFLY-PT} \rightarrow \textbf{knot type}$
\end{itemize}

The program \texttt{coords2egc} converts a single file of coordinates
to a single EGC.  We have included a python script
\texttt{batch\_coords2egc.py} on \cite{knottingtools} which can
be used to convert multiple coordinate files to a file of EGCs
(see Section \ref{sec:helperscripts} for more details).

The tools \texttt{xinger}, \texttt{jhomfly}, \texttt{jidknot},
\texttt{jidknot\_table16.py}, and \texttt{jidknot\_grep.py} all work
at the batch level, performing their duties on files with information
for one knot per line.  The file \texttt{jhomflytable.txt}
is needed by \texttt{jidknot\_table16.py} and \texttt{jidknot\_grep.py},
and is described in Section \ref{sec:prefixes}.

See Section \ref{sec:collisions} about
further steps that can be performed when there are collisions.

We describe how to download and install the software in the following
section.  In subsequent sections, we describe the steps of this
process from the end to the start (because it makes the explanation a
little easier).  Examples with full usage are in Section
\ref{sec:example}.

\subsection{Downloading and installing the software}
\label{sec:install}

Note that there are more thorough directions for downloading
and installing the software on \cite{knottingtools}.

First, the user should install \texttt{libhomfly}.  The
\texttt{libhomfly} package can be installed from repositories for all
of the major Linux distributions, or from GitHub \cite{libhomfly}.
For example, in Fedora Linux, typing the command
\begin{verbatim}
  $ sudo dnf install libhomfly
\end{verbatim}
will install the \texttt{libhomfly} package (as long as the user
has the appropriate permissions).

Second, the user should first visit \cite{knottingtools} to download the
packages for \texttt{coords2egc}, \texttt{xinger}, \texttt{jhomfly},
and \texttt{jidknot}.  These packages are compressed into \texttt{.tar.gz}
files.

Third, the user needs to install the packages for \texttt{coords2egc},
\texttt{xinger}, \texttt{jhomfly}, and \texttt{jidknot}.  For example,
after downloading a
package, say \texttt{coords2egc-somedate.tar.gz}, type
\begin{verbatim}
  $ tar -xf coords2egc-somedate.tar.gz
\end{verbatim}
to extract the package.  This will create a directory \texttt{COORDS2EGC},
which holds the code.  Then type
\begin{verbatim}
  $ cd COORDS2EGC/
  $ make
  $ mv coord2egc ~/bin/
  $ cd ../
\end{verbatim}
to make the executable \texttt{coords2egc} and then move that file to
one's local \texttt{bin/} directory.  Installation in Linux systems is
usually seamless, but not always.  We have included some installation
notes on \cite{knottingtools} with issues we have observed.  If a user
is having trouble with installation, please feel free to email either
of the authors.  Additional details can be found on the web page
\cite{knottingtools}.

Fourth, the user should download \texttt{jhomflytable.txt.gz} from
\cite{knottingtools} and extract the file \texttt{jhomflytable.txt}
using the commands
\begin{verbatim}
  $ gunzip jhomflytable.txt.gz
  $ mv jhomflytable.txt ~/bin/
\end{verbatim}
to extract the file and then move that file to one's local
\texttt{\$HOME/bin/} directory.  The user may also want to download the
alternatives to \texttt{jidknot} named \texttt{jidknot\_table16.py}
and/or \texttt{jidknot\_grep.py}.

After these steps, the user will have \texttt{libhomfly} installed on
their machine as well as working versions of \texttt{coords2egc},
\texttt{xinger}, \texttt{jhomfly}, and \texttt{jidknot}, each
in their local \texttt{\$HOME/bin/} directory.
The examples below assume that these programs are all
in the \texttt{\$HOME/bin/} directory and that the directory is within
their installation's local PATH variable.  Users can type
\begin{verbatim}
  $ echo $PATH
\end{verbatim}
to see those directories, and resources can be found on the web for
adding to the PATH.  Some helper scripts are described in Section
\ref{sec:helperscripts}.  Those scripts assume that these programs,
as well as \texttt{jhomflytable.txt}, are in the directory
\texttt{\$HOME/bin/}, although users can change this directory
in the scripts, if desired.

\subsection{Converting HOMFLY-PT polynomials to a knot type list using versions of \texttt{jidknot}}
\label{sec:jidknot}

There are three programs for converting a file of HOMFLY-PT polynomials
to a list of corresponding knot types: \texttt{jidknot},
\texttt{jidknot\_table16.py}, and \texttt{jidknot\_grep.py}.

The program \texttt{jidknot} is a simple C++ program that reads a file
containing plain-text versions of the HOMFLY-PT polynomial in the
LM-normalization output by \texttt{libhomfly}.  It uses an
\texttt{unordered\_map} from the C++ Standard Library to convert a
file of HOMFLY-PT polynomials to a list of chiral knot types matching
the given polynomials.  That container uses a hash table to provide
constant average time lookup.  The program \texttt{jidknot} has the
knot types through 12 crossings, and the links in the Rolfsen
table \cite{Rolfsen}.  Note that
\texttt{jidknot} takes a long time to compile.

A user who
needs access to knot types with crossing numbers through 16
should use one of the alternatives: \texttt{jidknot\_table16.py} or
\texttt{jidknot\_grep.py}.

The Python script \texttt{jidknot\_table16.py} converts a file of
HOMFLY-PT polynomials to knot/link types using a table file.  By
default, the program uses \texttt{jhomflytable.txt},
which includes all knot types through crossing number 16.  However, a
user could add or subtract from the table file, or create a new table
file.  The program loads the table file and creates a Python
dictionary for quick lookup.  Because the program creates the dictionary
on the fly, there is some startup cost (on EJR's current
machine using \texttt{jhomflytable.txt}, the startup takes around 2.5
seconds).  As such, this program is most useful for files with many
HOMFLY-PT polynomials.

If a user just has just a few HOMFLY-PT polynomials (or just one),
then \texttt{jidknot\_grep.py} might be a better choice.  The Python
script \texttt{jidknot\_grep.py} also converts a file of HOMFLY-PT
polynomials to knot/link types using a table file.  By default, again,
the program uses \texttt{jhomflytable.txt}.  This program utilizes the
Linux command \texttt{grep} to quickly scan the table file for
matches.

\subsection{Computing HOMFLY-PT polynomials using Jenkins's algorithm from crossing codes using \texttt{jhomfly}}

In his Master's thesis, Robert J.~Jenkins \cite{Jenkins} implemented a
load-balancing strategy to make the computation of the HOMFLY-PT
polynomial quicker, although the original implementation had some
problems that made it unreliable.  A repaired version of the code was
written by Jenkins in 2010 and was modified by Marco Miguel to be a
shared library \texttt{libhomfly}.\footnote{ Note that
  \texttt{libhomfly} is available from many Linux repositories.}
Miguel made some other changes to improve the usability of the
software and released it on GitHub \cite{libhomfly}, keeping the
public domain status.  The program \texttt{jhomfly}
\cite{knottingtools} uses \texttt{libhomfly} and reads extended Gauss
codes (EGCs) in the format described in \cite{ghomfly} (an example of
which is in Section \ref{sec:egc}), converting them to Jenkins's less
human-readable codes before sending them to \texttt{libhomfly}.  Then
\texttt{jhomfly} outputs the corresponding HOMFLY-PT polynomials in
the LM-normalization.  Note that \texttt{jhomfly} also fixes one
remaining issue with Jenkins's code in that it did not correctly
handle split diagrams (i.e., links containing a component that does
not cross any other component).

\subsection{Simplifying the crossing codes using \texttt{xinger}}

Since the HOMFLY-PT calculation is exponential\footnote{ This is only
  approximately true.  The computational complexity of the calculation
  is a ``complex'' question.  
  Jenkins gives the running time of his implementation as
  $O((m!)(2^n)(c^3))$ with memory space usage of $O((m!)(c^2))$ where
  $c$ is the number of crossings in the analyzed projection and $m$ is
  a quantity that has an upper bound of $\sqrt{c} + 1$
  \cite{Jenkins,JenkinsWebPage}.} on the number of crossings, it can
be useful to use the simplification schemes in \texttt{xinger} (based
on Reidemeister I and II moves) to remove extraneous crossings.  If
one is trying to compute knot types of random knots, for example, this
simplification step can save a significant amount of computation time.
If the configurations are fairly inane, then one can try to skip this
simplification step.  As such, this step is optional.  Also note that
\texttt{SnapPy} \cite{snappy} can perform simplification of crossing
codes (using Reidemeister I, II, and III moves) from planar diagram
(PD) codes.  However, a user would need to convert EGC to and from PD
codes to utilize this functionality from \texttt{SnapPy}.

\subsection{Converting 3D knot coordinates to crossing codes
  using \texttt{coords2egc}}
\label{sec:coords2egc}

\def\zoik#1{\texttt{\hspace*{2pc}#1}\\}

The program \texttt{coords2egc} can be used to create an extended
Gauss code (EGC, see Section \ref{sec:egc} below for a description)
from the vertex coordinates for a polygonal knot.  The coordinate file
can be in several different plain text formats including the Geomview
VECT format \cite{VECT} as well as the most basic format of the $xyz$
coordinates of one vertex per line (with a blank line separating
different components in the case of links).

For example, \\
\zoik{0.2 1.4 2.5}
\zoik{2.0 -1.0 -2.4}
\zoik{-3.5 -0.1 0.8}
\zoik{1.3 2.5 -2.2}
\zoik{-0.1 -3.9 0.5}
\zoik{0.2 3.6 1.0}
\zoik{2.6 -1.9 1.7}
\zoik{-2.7 -0.6 -2.0}
is a minimal stick representative of the $6_3$ knot type.

An example of a link file is\\
\zoik{2.9 -21.6 -9.0} 
\zoik{1.9 23.1 4.5} 
\zoik{-15.7 5.1 -8.2} 
\zoik{9.3 -7.1 14.3} 
\zoik{}
\zoik{14.5 -11.2 -4.7} 
\zoik{-13.9 -4.3 13.1} 
\zoik{1.1 16.0 -10.0}
which represents the link type $4^2_1$.



Note that the crossing code calculation in \texttt{coords2egc} is $O(n^2)$,
where $n$ is the number of edges in the polygon.  In theory, this
computation can be done in $O(n\log(n))$ time.  Currently, we are not
aware of any $O(n\log(n))$ implementations, although Rawdon's research
group is currently working on this problem.  If one is considering
long polygons with many crossings, it could be worth the effort to
implement a faster algorithm for translating vertex data to crossing
codes (although the HOMFLY-PT computation is likely to be the
bottleneck in large-scale knot type calculations).

We also have a Python script \texttt{batch\_coords2egc.py} which
converts multiple coordinate files to a file of EGCs,
which is described in Section \ref{sec:helperscripts}.

\subsection{Extended Gauss codes (EGCs)}
\label{sec:egc}

\begin{figure}[htb] 
\hbox to \textwidth {\hss\includegraphics[height=8cm]{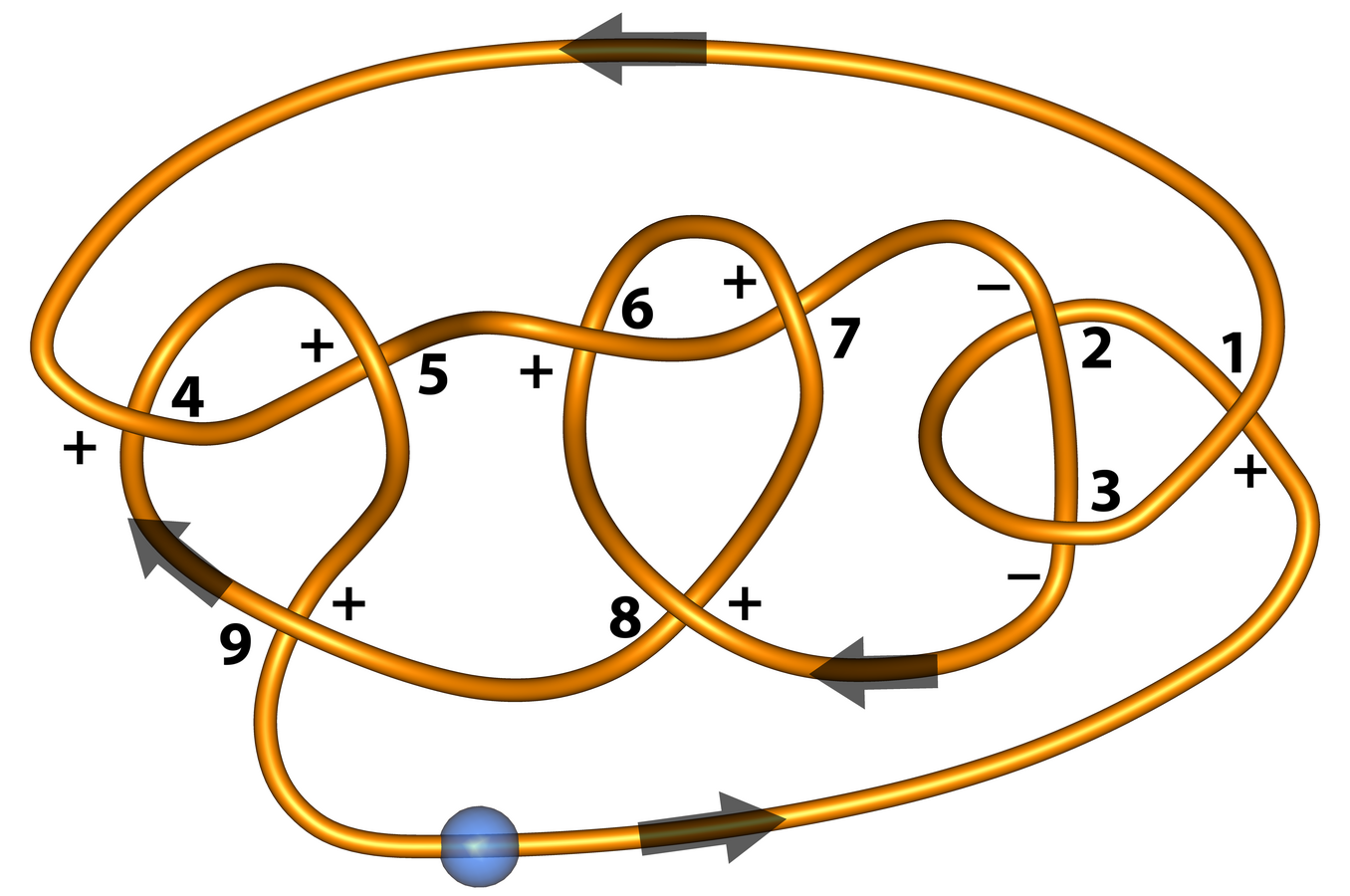}\hss}
\caption{Diagram of a $8_{21}$ knot with nine crossings labeled to produce
  the extended Gauss code (EGC) described in Section \ref{sec:example}.}
\label{fig:diagram}
\end{figure}

It may be the case that a user has their knots in some form other
than 3D coordinates.  In such a case, it may be easier for
the user to write a program that creates a file of crossing codes
directly.  Our programs work with extended Gauss codes (EGCs): the program \texttt{coords2egc} outputs an EGC
and \texttt{xinger} and \texttt{jhomfly} work from files of EGCs.

Here is a description of how to compute an EGC from the knot diagram in
Figure \ref{fig:diagram}.  First pick a starting point (the blue
sphere) and an orientation (denoted by the arrows).  Proceed along the
knot in the direction of the arrow and label each crossing with a
number in the order it occurs.  Also label each crossing according to
its sign.  Then start again and read off the crossing codes in
sequence, specifying above (\texttt{``a''}) if the strand is an over
crossing or below (\texttt{``b''}) if it goes under, the numerical
label, and the sign of the crossing.  For this example we get 
$$\texttt{b1+b2-a3-a1+a4+b5+a6+b7+a2-b3-a8+b6+a7+b8+a9+b4+a5+b9+}$$
as the EGC.

\section{Examples}
\label{sec:example}

We provide examples for using each individual program as well as
all of them together.

We assume below that the executables are
all within the user's path, ideally in \texttt{\$HOME/bin/} as mentioned
in Section \ref{sec:install}.  If the executables are not in one's
path, substitute
\begin{verbatim}
  /path/to/coords2egc
\end{verbatim}
for
\begin{verbatim}
  coords2egc
\end{verbatim}
where \texttt{/path/to/} is the directory holding the given program.

\subsection{Convert a knot/link file into an EGC using \texttt{coords2egc}}
\label{sec:coords2egcexample}

First, save the coordinates for the $6_3$ knot in Section \ref{sec:coords2egc}
to a file named \texttt{6.3.txt}.

The command
\begin{verbatim}
  $ coords2egc 6.3.txt
\end{verbatim}
will output
$$\texttt{a1-b2-b3+a4-b5-b6+b7-b1-a2-b8+a9+a3+b10+a5-a6+a7-a8+b9+b4-a10+}$$
to the terminal.  Use
\begin{verbatim}
  $ coords2egc 6.3.txt 6.3.egc
\end{verbatim}
to save the output to a file named \texttt{6.3.egc}.

\subsection{Simplify a file of EGC codes using \texttt{xinger}}
\label{sec:xingerexample}

Save the following lines to a file named \texttt{egcs.txt}.
\begin{verbatim}
b1-b2+b3-a3-a2+a1-
a1-b2-b3+a4-b5-b6+b7-b1-a2-b8+a9+a3+b10+a5-a6+a7-a8+b9+b4-a10+
\end{verbatim}

The command
\begin{verbatim}
  $ xinger +r egcs.txt
\end{verbatim}
will output
\begin{verbatim}
b1-a1-
a1-b2-b3+a4-b5-b1-a2-b6+a7+a3+b8+a5-a6+b7+b4-a8+
\end{verbatim}
to the terminal.  Note that the first of these EGCs is for an unknot, so
the returned EGC forms a simple twist.  These twists are kept so that
we do not lose track of components, which is critical for unknots and
split links.

As above, one can output to a file, here \texttt{simplified\_egcs.txt},
using the following.
\begin{verbatim}
  $ xinger +r egcs.txt simplified_egcs.txt
\end{verbatim}

\subsection{Compute HOMFLY-PT polynomials from a file of EGCs using \texttt{jhomfly}}
\label{sec:homflyptstxt}

Using the file \texttt{simplified\_egcs.txt} from the previous section,
the command
\begin{verbatim}
  $ jhomfly simplified_egcs.txt homflypts.txt
\end{verbatim}
will place
\begin{verbatim}
1
L^-2 + 3 + L^2 - M^2L^-2 - 3M^2 - M^2L^2 + M^4
\end{verbatim}
in the file \texttt{homflypts.txt}.  Similar to above, the user
can use \texttt{jhomfly} without a second argument to print the results
to the terminal.

\subsection{Translate a file of HOMFLY-PT polynomials to a file of knot types}

Recall that there are three options, which are described in Section \ref{sec:jidknot}.

Using the file \texttt{homflypts.txt} from the previous section,
the commands
\begin{verbatim}
  $ jidknot -k homflypts.txt knot_types.txt
  $ jidknot_table16.py homflypts.txt > knot_types.txt
  $ jidknot_grep.py homflypts.txt > knot_types.txt
\end{verbatim}
will place
\begin{verbatim}
a0.1
a6.3
\end{verbatim}
in the file \texttt{knot\_types.txt}.  The prefixes of ``\texttt{a}''
denote that the knot types are achiral.  The other prefixes are
``\texttt{p}'', ``\texttt{m}'', and ``\texttt{h}'', and are explained
in Section \ref{sec:prefixes}.  For the two Python scripts, the user
must change permissions to make the scripts executable, or else
prepend the call with \texttt{python} as follows.
\begin{verbatim}
  $ python jidknot_table16.py homflypts.txt > knot_types.txt
\end{verbatim}

\subsection{Using piping to string together the programs}

Users can use POSIX pipes to avoid intermediate files and go straight
from a coordinate file to a knot type.  For all of the programs mentioned
above, if the program receives one argument, the output is sent to
STDOUT, i.e., the terminal.  If a second argument is given, then
the output is sent to the file name given by the second argument.
The \texttt{-{}-} below tell the programs to receive input from
STDIN instead of a file.

Using \texttt{6.3.txt}
from Section \ref{sec:coords2egcexample}, the command
\begin{verbatim}
  $ coords2egc 6.3.txt | xinger +r -- | jhomfly -- | jidknot -k --
\end{verbatim}
outputs \texttt{a6.3}.  Users can remove the ``\texttt{| xinger +r -{}-}''
portion to skip the simplification step.

Working from the file \texttt{egcs.txt} from Section \ref{sec:xingerexample},
the commands
\begin{verbatim}
  $ xinger +r egcs.txt | jhomfly -- | jidknot -k --
\end{verbatim}
or
\begin{verbatim}
  $ jhomfly egcs.txt | jidknot -k --
\end{verbatim}
both output
\begin{verbatim}
a0.1
a6.3
\end{verbatim}
which is the list of knot types.

\subsection{Helper scripts}
\label{sec:helperscripts}

We have included four helper scripts written in Python which
just string together the commands using pipes (like in the example in
the previous section).  More information can be found on the site
\cite{knottingtools}.

Each of these scripts is written in Python and is self contained.
One can download the scripts of interest, or none at all.

\subsubsection{Technical notes}

These scripts each require that \texttt{libhomfly} is installed and
that the executables 
\texttt{coords2egc}, \texttt{jhomfly}, and \texttt{jidknot} are
in \texttt{\$HOME/bin/}.  One can change the directory for the
executables in the scripts if need be.  In some cases, not all of the programs
are needed.

For some of the scripts, the file \texttt{jhomflytable.txt} is
required and that the file is also in \texttt{\$HOME/bin/}.

We also assume that the user has placed these scripts in
\texttt{\$HOME/bin/}, that \texttt{\$HOME/bin/} is in the user's
\texttt{PATH}, and that the user has changed the permissions of the
scripts to be executable.  Otherwise, the user must prepend
\texttt{python} to all of the script calls, for example
\begin{verbatim}
  $ python coords2knottype.py coordinates-file
\end{verbatim}
to match the call in the next subsection.

There is more information about these scripts and examples of their
usage at \cite{knottingtools}.

\subsubsection{Convert a single coordinate file to a knot type}
\begin{verbatim}
  $ coords2knottype.py [-f other_table_file] [-j] coordinates-file
\end{verbatim}
Here \texttt{coordinates-file} is a file of vertex coordinates in a
supported format, e.g., the plain text file mentioned in Section
\ref{sec:coords2egc}.  This script uses the code base from
\texttt{jidknot\_grep.py} and the table file \texttt{jhomflytable.txt}
for the knot type lookup (all knot types through crossing number 16).
The optional flag \texttt{-f other\_table\_file} loads a user specified
table file in place of \texttt{jhomflytable.txt}.
The user can use the optional \texttt{-j} flag to force the use of
\texttt{jidknot}, which does not require a table file and
contains knot types (and some link types) through crossing number 12.

\subsubsection{Convert a number of coordinate files to EGCs}
\begin{verbatim}
  $ batch_coords2egc.py file-of-file-names > file-of-EGCs
\end{verbatim}
For this script, the file \texttt{file-of-file-names} contains a list of
coordinate file names, one per line.  This script writes the EGCs to
STDOUT, which, in this example, is then piped to the file
\texttt{file-of-EGCs}.

\subsubsection{Convert a number of coordinate files to knot types} 
\begin{verbatim}
  $ batch_coords2knottypes.py [-n] [-j|-f other_table_file]
      file-of-file-names > file-of-knot-types
\end{verbatim}
Again, the file \texttt{file-of-file-names} contains a list of
coordinate file names, one per line.  This script writes the knot
types to STDOUT, which, in this example, is then piped to the file
\texttt{file-of-knot-types}.  By default, this program simplifies the
EGCs using \texttt{xinger} and utilizes the code base from
\texttt{jidknot\_table16.py} (which requires a table file).  The
optional flag \texttt{-n} flag skips the \texttt{xinger}
simplification.  The table file defaults to \texttt{jhomflytable.txt}
but a different table file can be specified with the optional flag
\texttt{-f other\_table\_file}.  The optional flag \texttt{-j} forces
the use of \texttt{jidknot}.

\subsubsection{Convert a file of EGCs to knot types}
\begin{verbatim}
  $ egc2knottypes.py [-n] [-j|-f other_table_file]
      file-of-EGCs > file-of-knot-types
\end{verbatim}
The program has the same default choices as
\texttt{batch\_coords2knottypes.py} in the previous section and also
supports the optional flags \texttt{-n}, \texttt{-j}, and
\texttt{-f}.

\section{Information about knot types through 16 crossings}

In this section, we provide some general information about knot
types through crossing number 16.  None of this information is new, 
but nonetheless could be useful to the reader.


\subsection{Achiral prime knot types}
Through 10 crossings, there are 21 prime knot types which are achiral:
$0_1$, 
$4_1$, 
$6_3$, 
$8_3$, 
$8_9$, 
$8_{12}$, 
$8_{17}$, 
$8_{18}$, 
$10_{17}$, 
$10_{33}$, 
$10_{37}$, 
$10_{43}$, 
$10_{45}$, 
$10_{79}$, 
$10_{81}$, 
$10_{88}$, 
$10_{99}$, 
$10_{109}$, 
$10_{115}$, 
$10_{118}$, and
$10_{123}$.
From \cite{HosThiWee}, we know that there are
0, 58, 0, 274, 1, and 1539
prime achiral knot types with crossing numbers
11 through 16, respectively.\footnote{We thank Morwen Thistlethwaite
  for providing the list of achiral knot types from 11 through
  16 crossings.}


\subsection{Chiral knot types whose mirror pair share the same HOMFLY-PT polynomial}
Through 10 crossings, there are six prime knot types for which the base
knot type is chiral but the HOMFLY-PT polynomial
cannot distinguish between the two chiralities, namely:
$9_{42}$, 
$10_{48}$, 
$10_{71}$, 
$10_{91}$, 
$10_{104}$, and
$10_{125}$.
There are 2, 33, 35, 342, 394, and 2660 prime knot types with crossing numbers
11 through 16, respectively, with this property.

\subsection{Number of HOMFLY-PT polynomials and collisions}

Through 16 crossings, there are 2,427,782 unique HOMFLY-PT polynomials:
2,367,155 containing only prime knot types, 58,237 containing only composite
knot types, and 2390 containing a mixture of prime and composite knot types.

In the introduction, we say that the HOMFLY-PT polynomial provides a
good balance of speed and effectiveness in identifying knot types
through 16 crossings.  Of course, this depends on one's
tolerance for collisions.  In Table \ref{collisiontable2}, we show the
efficacy of HOMFLY-PT in identifying knot types (both prime and
composite) through different crossing numbers.  The table shows that
if one were to know, for example, that all knot types observed have
crossing number eight or smaller, then HOMFLY-PT correctly identifies
all chiral knot types.  Through crossing number 10, there are only 18
collisions, and each of these 18 collisions only contain a single pair of two chiral
knot types.  One can see that through crossing number 16 knot types,
approximately 74\% of the observed HOMFLY-PT polynomials have unique
chiral knot types associated with them.  However, at crossing number 16, we
also see that there is one HOMFLY-PT polynomial containing 40 matching
chiral knot types.  Thus, one can use Table \ref{collisiontable2} to
get some sense for the number of HOMFLY-PT polynomial collisions over
different crossing numbers.


\begin{table}
  \caption{The number of HOMFLY-PT polynomials with different numbers
    of collisions (or non-collisions) through the given crossing
    number (xing).  NCKT is the number of chiral knot types through
    the given crossing number.  NHF is the number of HOMFLY-PT polynomials observed
    through the given crossing number. The rows below show the breakdown
    of how many HOMFLY-PT polynomials observed through the given
    crossing number have that many chiral knot types in them.
    For example, for HOMFLY-PT polynomials observed with crossing number
    14 or smaller, there are 2957 HOMFLY-PT polynomials which have exactly three
    chiral knot types sharing a HOMFLY-PT polynomial.}
  \label{collisiontable2}
  \centering
  \begin{tabular}{c||c|c|c|c|c|c|c|c|c}
    xing &  8 &   9 &  10 &   11 &   12 &     13 &     14 &      15 &        16 \\
    \hline
    NCKT & 78 & 194 & 557 & 1783 & 6410 & 27,416 & 124,506& 643,663 & 3,473,712 \\
    \hline
    NHF  & 78 & 191 & 539 & 1669 & 5717 & 22,498 & 95,904 & 464,222 & 2,427,782 \\
    \hline
    \hline
    1    & 78 & 188 & 521 & 1559 & 5094 & 18,644 & 75,220 & 347,788 & 1,787,634 \\
    2    &    & 3   & 18  & 106  & 560  & 3101   & 15,866 & 83,517  & 450,557   \\
    3    &    &     &     & 4    & 58   & 534    & 2957   & 16,807  & 82,463    \\
    4    &    &     &     &      & 4    & 157    & 1206   & 9915    & 66,435    \\
    5    &    &     &     &      &      & 40     & 351    & 2595    & 14,881    \\
    6    &    &     &     &      & 1    & 17     & 179    & 1847    & 11,480    \\
    7    &    &     &     &      &      & 4      & 42     & 677     & 4082      \\
    8    &    &     &     &      &      &        & 51     & 586     & 5049      \\
    9    &    &     &     &      &      &        & 14     & 160     & 1413      \\
    10   &    &     &     &      &      & 1      & 10     & 155     & 1285      \\
    11   &    &     &     &      &      &        & 6      & 48      & 574       \\
    12   &    &     &     &      &      &        & 12     & 60      & 771       \\
    13   &    &     &     &      &      &        &        & 20      & 261       \\
    14   &    &     &     &      &      &        &        & 14      & 264       \\
    15   &    &     &     &      &      &        &        & 4       & 132       \\
    16   &    &     &     &      &      &        &        & 15      & 157       \\
    17   &    &     &     &      &      &        &        & 6       & 69        \\
    18   &    &     &     &      &      &        &        & 2       & 75        \\
    19   &    &     &     &      &      &        &        &         & 47        \\
    20   &    &     &     &      &      &        &        &         & 49        \\
    21   &    &     &     &      &      &        &        &         & 24        \\
    22   &    &     &     &      &      &        &        & 1       & 21        \\
    23   &    &     &     &      &      &        &        & 4       & 18        \\
    24   &    &     &     &      &      &        & 1      &         & 17        \\
    25   &    &     &     &      &      &        &        &         & 3         \\
    26   &    &     &     &      &      &        &        &         & 8         \\
    27   &    &     &     &      &      &        &        &         & 2         \\
    28   &    &     &     &      &      &        &        &         & 4         \\
    29   &    &     &     &      &      &        &        &         & 1         \\
    30   &    &     &     &      &      &        &        & 1       & 1         \\
    32   &    &     &     &      &      &        &        &         & 1         \\
    35   &    &     &     &      &      &        &        &         & 2         \\
    36   &    &     &     &      &      &        &        &         & 1         \\
    40   &    &     &     &      &      &        &        &         & 1         \\
  \end{tabular}
  \label{collisioncounts2}
\end{table}

\subsection{Peculiarities in the knot tables}
\label{sec:peculiarities}

In this section, we note some potential discrepancies between the
naming of knot types when using different pieces of software.

First, there is some disagreement in software and books between the
$10_{83}$ and $10_{86}$ knot types.  In the original edition of
Rolfsen's Knots and Links book \cite{Rolfsen}, the Conway notation and
Alexander polynomial listed for $10_{83}$ is really associated with
the picture for $10_{86}$, and vice versa.  Depending on the software
one is using, there will be different interpretations of which knot type
is $10_{83}$ and which is $10_{86}$.

Second, Rolfsen's original table had 166 knot types listed under
crossing number 10.  Perko \cite{perko} discovered that two of the
knot types were actually equivalent, the (\textbf{former})
$\mathbf{10_{161}}$ and $\mathbf{10_{162}}$.  These two types were
joined and are now called $10_{161}$.  The knot types formerly known as
$\mathbf{10_{163}}$ through $\mathbf{10_{166}}$ had each of their
indices decreased by one to form $10_{162}$ through $10_{165}$.  Note
that, again, we see different behavior from different programs.  For
example, \texttt{SnapPy} shows all 166 whereas KnotPlot uses
the new notation.

Third, through crossing number 10, researchers typically use the
so-called Alexander-Briggs-Conway (ABC) notation for prime knot types,
which is simply $xing_{index}$, where $xing$ is the crossing number
and $index$ is the traditional placement in the knot table amongst
knot types with the given crossing number.  For crossing numbers 11
through 16 \cite{HosThiWee}, the Dowker-Thistlethwaite (DT) notation
is used. Here the index is based on lexicographical ordering of the
minimal DT notation, and alternating and nonalternating knot types are
enumerated separately.  For example, this notation looks like $11a99$
for the $99^{\text{th}}$ alternating knot type with crossing number 11
and $11n27$ for the $27^{\text{th}}$ nonalternating knot type with
crossing number 11.  The knot types through crossing number 10 also
have DT notation, although the indices differ from the ABC notation.
Note that in enumerating the prime knot types through crossing number
19, Burton \cite{Burton} has a different notation, separating the knot
types into the classes: alternating hyperbolic, nonalternating
hyperbolic, torus, and satellite knot types.

\subsection{Counts for chiral knot types}
\label{sec:chiralknottypes}

Since the HOMFLY-PT polynomials usually provides data about chiral
knot types yet the traditional knot tables enumerate base knot types
(i.e., are agnostic with respect to
mirror images), we provide some tables of chiral knot type counts by
crossing number.  Table \ref{primetable} shows the counts for prime
knot types through crossing number 16.  Tables \ref{comp2table},
\ref{comp3table}, \ref{comp4table}, and \ref{comp5table} show the
number of chiral composite knot types by the sum of the factor knots
(which we assume to be the crossing number).

\begin{table}
  \caption{Number of prime chiral knot types by crossing number
    divided into the categories of chiral knot type pairs with
    distinct HOMFLY-PT polynomials, chiral knot type pairs whose pairs
    share a common HOMFLY-PT polynomial, and achiral knot types.  Each
    of these categories is then subdivided into the number of
    alternating (alt) and nonalternating (non) knot types with the
    given property.  The final column shows the total number of chiral
    knot types, which is the sum of columns 6 and 7 plus twice the sum
    of columns 2 through 5.  The sum of columns 2 through 7 gives the
    traditional counts of base knot types from the knot tables of
    \cite{HosThiWee}.}
  \label{primetable}
  \vspace{0.1in}
  \centering
  \begin{tabular}{c||c|c|c|c|c|c||c}
    & \multicolumn{4}{|c|}{chiral pairs} & \multicolumn{2}{|c||}{achiral types} & chiral\\
    \hline
    & \multicolumn{2}{|c|}{distinct HF} & \multicolumn{2}{|c|}{common HF} & \multicolumn{2}{|c||}{} & knot types\\
    \hline
    crossing & alt     & non       & alt  & non  & alt  & non &           \\
    \hline                                                        
    0        &         &           &      &      & 1    &     & 1         \\
    3        & 1       &           &      &      &      &     & 2         \\
    4        &         &           &      &      & 1    &     & 1         \\
    5        & 2       &           &      &      &      &     & 4         \\ 
    6        & 2       &           &      &      & 1    &     & 5         \\ 
    7        & 7       &           &      &      &      &     & 14        \\ 
    8        & 13      & 3         &      &      & 5    &     & 37        \\ 
    9        & 41      & 7         &      & 1    &      &     & 98        \\ 
    10       & 106     & 41        & 4    & 1    & 13   &     & 317       \\ 
    11       & 367     & 183       &      & 2    &      &     & 1104      \\ 
    12       & 1212    & 873       & 22   & 11   & 54   & 4   & 4294      \\ 
    13       & 4878    & 5075      &      & 35   &      &     & 19,976    \\ 
    14       & 19,140  & 27,216    & 173  & 169  & 223  & 51  & 93,670    \\ 
    15       & 85,263  & 167,635   &      & 394  &      & 1   & 506,585   \\ 
    16       & 377,708 & 1,006,798 & 1042 & 1618 & 1049 & 490 & 2,775,871 \\
    \hline                                                        
    total    & 488,740 & 1,207,831 & 1241 & 2231 & 1347 & 546 & 3,401,979 \\ 
  \end{tabular}
\end{table}

\begin{table}
  \caption{Counts for number of composite chiral knot types with two
    factors by crossing number.  Bold indicates that the total
    includes some compositions containing chiral knot types whose
    mirror images share the same HOMFLY-PT polynomial.  The second and
    third columns show the crossing numbers of the prime factor knots.
    For example, the line containing ``11 5 6 20'' says that there are
    20 chiral composite knot types with crossing number 11 formed by
    composing 5- and 6-crossing prime knot types.}
  \label{comp2table}
  \vspace{0.1in}
  \centering
  \begin{tabular}{c||c|c||c||c}
    & \multicolumn{2}{|c||}{crossing \#} & \multicolumn{2}{c}{} \\
    & \multicolumn{2}{|c||}{of factor}   & \multicolumn{2}{c}{} \\
    \hline
    crossings & 1   & 2   & count           & total           \\
    \hline
    \hline
    6         & {\ }{\ }{\ }3{\ }{\ }{\ }   & {\ }{\ }{\ }3{\ }{\ }{\ }   & 3               & 3               \\
    \hline
    7         & 3   & 4   & 2               & 2               \\
    \hline
    8         & 3   & 5   & 8               & 9               \\
    8         & 4   & 4   & 1               &                 \\
    \hline
    9         & 3   & 6   & 10              & 14              \\
    9         & 4   & 5   & 4               &                 \\
    \hline
    10        & 3   & 7   & 28              & 43              \\
    10        & 4   & 6   & 5               &                 \\
    10        & 5   & 5   & 10              &                 \\
    \hline
    11        & 3   & 8   & 74              & 108             \\
    11        & 4   & 7   & 14              &                 \\
    11        & 5   & 6   & 20              &                 \\
    \hline
    12        & 3   & 9   & \textbf{196}    & \textbf{304}    \\
    12        & 4   & 8   & 37              &                 \\
    12        & 5   & 7   & 56              &                 \\
    12        & 6   & 6   & 15              &                 \\
    \hline
    13        & 3   & 10  & \textbf{634}    & \textbf{950}    \\
    13        & 4   & 9   & \textbf{98}     &                 \\
    13        & 5   & 8   & 148             &                 \\
    13        & 6   & 7   & 70              &                 \\
    \hline
    14        & 3   & 11  & \textbf{2208}   & \textbf{3207}   \\
    14        & 4   & 10  & \textbf{317}    &                 \\
    14        & 5   & 9   & \textbf{392}    &                 \\
    14        & 6   & 8   & 70              &                 \\
    14        & 7   & 7   & 105             &                 \\
    \hline
    15        & 3   & 12  & \textbf{8588}   & \textbf{11,968} \\
    15        & 4   & 11  & \textbf{1104}   &                 \\
    15        & 5   & 10  & \textbf{1268}   &                 \\
    15        & 6   & 9   & \textbf{490}    &                 \\
    15        & 7   & 8   & 518             &                 \\
    \hline
    16        & 3   & 13  & {\ }\textbf{39,952}{\ } & \textbf{52,322} \\
    16        & 4   & 12  & \textbf{4294}   &                 \\
    16        & 5   & 11  & \textbf{4416}   &                 \\
    16        & 6   & 10  & \textbf{1585}   &                 \\
    16        & 7   & 9   & \textbf{1372}   &                 \\
    16        & 8   & 8   & 703             &                 \\
    \hline
    total     &     &     &                 & {\ }\textbf{68,930}{\ } \\
  \end{tabular}
\end{table}

\begin{table}
  \caption{Counts for number of composite chiral knot types with three
    factors by crossing number.  Bold indicates that the total
    includes some compositions containing chiral knot types whose
    mirror images share the same HOMFLY-PT polynomial.  The second,
    third, and fourth columns show the crossing numbers of the prime
    factor knots.}
  \label{comp3table}
  \vspace{0.1in}
  \centering
  \begin{tabular}{c||c|c|c||c||c}
    & \multicolumn{3}{|c||}{crossing \#} & \multicolumn{2}{c}{} \\
    & \multicolumn{3}{|c||}{of factor}   & \multicolumn{2}{c}{} \\
    \hline
    crossings & 1 & 2 & 3 & count        & total         \\
    \hline
    9         & {\ }{\ }3{\ }{\ }   & {\ }{\ }3{\ }{\ }   & {\ }{\ }3{\ }{\ }   & 4            & 4             \\
    \hline
    10        & 3   & 3   & 4   & 3            & 3             \\
    \hline
    11        & 3   & 3   & 5   & 12           & 14            \\
    11        & 3   & 4   & 4   & 2            &               \\
    \hline
    12        & 3   & 3   & 6   & 15           & 24            \\
    12        & 3   & 4   & 5   & 8            &               \\
    12        & 4   & 4   & 4   & 1            &               \\
    \hline
    13        & 3   & 3   & 7   & 42           & 76            \\
    13        & 3   & 4   & 6   & 10           &               \\
    13        & 3   & 5   & 5   & 20           &               \\
    13        & 4   & 4   & 5   & 4            &               \\
    \hline
    14        & 3   & 3   & 8   & 111          & 194           \\
    14        & 3   & 4   & 7   & 28           &               \\
    14        & 3   & 5   & 6   & 40           &               \\
    14        & 4   & 4   & 6   & 5            &               \\
    14        & 4   & 5   & 5   & 10           &               \\
    \hline
    15        & 3   & 3   & 9   & \textbf{294} & \textbf{564}  \\
    15        & 3   & 4   & 8   & 74           &               \\
    15        & 3   & 5   & 7   & 112          &               \\
    15        & 3   & 6   & 6   & 30           &               \\
    15        & 4   & 4   & 7   & 14           &               \\
    15        & 4   & 5   & 6   & 20           &               \\
    15        & 5   & 5   & 5   & 20           &               \\
    \hline
    16        & 3   & 3   & 10  & \textbf{951} & \textbf{1741} \\
    16        & 3   & 4   & 9   & \textbf{196} &               \\
    16        & 3   & 5   & 8   & 296          &               \\
    16        & 3   & 6   & 7   & 140          &               \\
    16        & 4   & 4   & 8   & 37           &               \\
    16        & 4   & 5   & 7   & 56           &               \\
    16        & 4   & 6   & 6   & 15           &               \\
    16        & 5   & 5   & 6   & 50           &               \\
    \hline
    total     &     &     &     &              & {\ }\textbf{2620}{\ } \\
  \end{tabular}
\end{table}

\begin{table}
  \caption{Counts for number of composite chiral knot types with four
    factors by crossing number. The second through fifth columns show
    the crossing numbers of the prime factor knots.}
  \label{comp4table}
  \vspace{0.1in}
  \centering
  \begin{tabular}{c||c|c|c|c||c||c}
    & \multicolumn{4}{|c||}{crossing \#} & \multicolumn{2}{c}{} \\
    & \multicolumn{4}{|c||}{of factor}   & \multicolumn{2}{c}{} \\
    \hline
    crossings & 1   & 2   & 3   & 4   & count        & total \\
    \hline
    \hline
    12        & {\ }3{\ }   & {\ }3{\ }   & {\ }3{\ }   & {\ }3{\ }   & 5            & 5     \\
    \hline
    13        & 3   & 3   & 3   & 4   & 4            & 4     \\
    \hline
    14        & 3   & 3   & 3   & 5   & 16           & 19    \\
    14        & 3   & 3   & 4   & 4   & 3            &       \\
    \hline
    15        & 3   & 3   & 3   & 6   & 20           & 34    \\
    15        & 3   & 3   & 4   & 5   & 12           &       \\
    15        & 3   & 4   & 4   & 4   & 2            &       \\
    \hline
    16        & 3   & 3   & 3   & 7   & 56           & 110   \\
    16        & 3   & 3   & 4   & 6   & 15           &       \\
    16        & 3   & 3   & 5   & 5   & 30           &       \\
    16        & 3   & 4   & 4   & 5   & 8            &       \\
    16        & 4   & 4   & 4   & 4   & 1            &       \\
    \hline
    total     &     &     &     &     &              & 172   \\
  \end{tabular}
\end{table}

\begin{table}
  \caption{Counts for number of composite chiral knot types with five factors by crossing number. The second through sixth columns show
    the crossing numbers of the prime factor knots.}
  \label{comp5table}
  \vspace{0.1in}
  \centering
  \begin{tabular}{c||c|c|c|c|c||c||c}
    & \multicolumn{5}{|c||}{crossing \#} & \multicolumn{2}{c}{} \\
    & \multicolumn{5}{|c||}{of factor}   & \multicolumn{2}{c}{} \\
    \hline
    crossings & 1 & 2 & 3 & 4 & 5 & count  & total \\
    \hline
    \hline
    15        & {\ }3{\ }   & {\ }3{\ }   & {\ }3{\ }   & {\ }3{\ }   & {\ }3{\ }   & 6      & 6     \\
    \hline
    16        & 3   & 3   & 3   & 3   & 4   & 5      & 5     \\
    \hline
    total     &     &     &     &     &     &        & 11    \\
   \end{tabular}
\end{table}
                                                      
\begin{table}
  \caption{Counts for number of chiral knot types by crossing number.}
  \label{final}
  \vspace{0.1in}
  \centering
  \begin{tabular}{c||c|c|c|c|c||c}
    \multicolumn{2}{c}{} & \multicolumn{4}{|c||}{number of composite factors} & \\
    \hline
    crossing & prime     &  2     &  3   &  4   & 5    & total       \\
    \hline            
    0        & 1         &        &      &      &      & 1           \\
    3        & 2         &        &      &      &      & 2           \\    
    4        & 1         &        &      &      &      & 1           \\  
    5        & 4         &        &      &      &      & 4           \\ 
    6        & 5         & 3      &      &      &      & 8           \\ 
    7        & 14        & 2      &      &      &      & 16          \\ 
    8        & 37        & 9      &      &      &      & 46          \\ 
    9        & 98        & 14     & 4    &      &      & 116         \\ 
    10       & 317       & 43     & 3    &      &      & 363         \\ 
    11       & 1104      & 108    & 14   &      &      & 1226        \\ 
    12       & 4294      & 304    & 24   & 5    &      & 4627        \\ 
    13       & 19,976    & 950    & 76   & 4    &      & 21,006      \\ 
    14       & 93,670    & 3207   & 194  & 19   &      & 97,090      \\ 
    15       & 506,585   & 11,968 & 564  & 34   & 6    & 519,157     \\ 
    16       & 2,775,871 & 52,322 & 1741 & 110  & 5    & 2,830,049   \\
    \hline            
    total    & {\ }3,401,979{\ } & {\ }68,930{\ } & {\ }2620{\ } & {\ }172{\ }  & {\ }11{\ }   & {\ }3,473,712{\ }   \\
  \end{tabular}
\end{table}

\section{Resolving some types of HOMFLY-PT collisions}
\label{sec:collisions}

Earlier we noted that HOMFLY-PT collisions can occur for even
small numbers of crossings, e.g., $+5_1$ and $+10_{132}$ share
the same HOMFLY-PT polynomial.  If one were to want to determine
whether a given configuration forms $+5_1$ versus $+10_{132}$,
more work is needed.  There are a number of techniques that
can be used to deal with such issues.  Below we present
some thoughts along these lines.  This section is not meant
to be complete, rather a list of some observations about what
we, and others, have tried to resolve these issues.

Imagine we have a configuration whose HOMFLY-PT polynomial matches
both $+5_1$ and $+10_{132}$.  The simplest situation would be if the
crossing code could be simplified (using \texttt{xinger} or some other
program) to 9 or fewer crossings, in which case the knot type must be
$+5_1$.  However, if the crossing code could be simplified only to 10 or
greater crossings, we cannot definitively conclude the knot type.  It could
be that the knot type is still $+5_1$, but the simplification algorithm
failed to reduce the crossing code to below 10 crossings.  Or the
configuration could form a knot type whose crossing number is 17 or
higher.

One strategy is to use different knot invariants to try to identify
the knot type.  There is a host of topological information about prime
knot types through 12 crossings at KnotInfo \cite{knotinfo}.  Also,
the software \texttt{pyknotid} \cite{pyknotid} has knot invariant
information for prime knot types with crossing number up to
15. Invariants such as the hyperbolic volume (when the knot types are
hyperbolic) can then be used to rule out (up to computational
inaccuracies) certain knot types.

A second, and more robust, strategy is to identify the base knot type.
There are algorithms that can identify the exact base knot type
for many configurations.  If one only cares about the base knot
types, then these algorithms are sufficient.  Otherwise, the
information can be used in conjunction with the HOMFLY-PT polynomial
to identify chiral knot types.

For example, the program \texttt{knotfind}, a part of the (now
unsupported) \texttt{Knotscape} package \cite{knotscape}, uses Dowker
codes to identify the exact base prime knot type for knot types through
crossing number 16.  In our experience, the full \texttt{Knotscape} package
does not fully compile on Fedora Linux, but the compiling works far
enough to generate \texttt{knotfind}.  The program works by
converting a Dowker code into a reduced Dowker code, which can then be
matched to data provided in the package to find the knot type.  If
\texttt{knotfind} finds that a Dowker code corresponds to a
composition, it factors the Dowker code into the Dowker codes for the
two factors.  Note that \texttt{knotfind} can be finicky at times, the
details of which are too involved to explain here.  If one wants to
pursue this option, please contact author Rawdon.

Elements of \texttt{SnapPy} \cite{snappy} can also be used to compute
exact base prime knot types.  In particular, author Rawdon has used
\texttt{SnapPy} to calculate isometry signatures, which can then be used
to determine exact base knot types.  The isometry signature is only defined
for hyperbolic knot types, and it is critical to note that composite
knot types are not hyperbolic.  Burton used canonical triangulation
calculations in \texttt{SnapPy} in enumerating the knot types with
crossing number between 17 and 19 \cite{Burton}, so this would be
another avenue one could explore.

We note one curious case that arose recently with respect
to using these base knot type techniques.  Recall that $+5_1$ and
$+10_{132}$ share the same HOMFLY-PT polynomial.  If we
have a knot configuration $K$ with this HOMFLY-PT polynomial
and \texttt{knotfind} returns that the base knot type of $K$ is
$10_{132}$, then we know that the exact chiral knot type is
$+10_{132}$.
But now consider the knot types $+5_1\#-10_{132}$ and
$-5_1\#+10_{132}$, both of whose base knot type is $5_1\#10_{132}$.
Due to the multiplicative nature of the HOMFLY-PT polynomial, the
polynomials for $+5_1\#-10_{132}$ and $-5_1\#+10_{132}$ also match.
Thus, the base knot type and HOMFLY-PT polynomial information together
is not sufficient to determine which of $+5_1\#-10_{132}$ or
$-5_1\#+10_{132}$ the configuration forms.  Furthermore, the knot
types for $+3.1\#-13n585$ and $-3.1\#+13n585$ match, which creates an
analogous problem.  Through crossing number 16, these two are the only
cases observed where (beyond the chiral knot types containing knot
types with matching HOMFLY-PT polynomials over mirror pairs) the
combination of base knot type and HOMFLY-PT polynomial is not
sufficient to determine the exact chiral knot type.

\section{Our \texttt{libhomfly} to knot types conversion file}
\label{sec:prefixes}

The file \texttt{jhomflytable.txt} at \cite{knottingtools} contains
the translation from \texttt{jhomfly} output (which comes from
\texttt{libhomfly}) to knot types matching
the HOMFLY-PT polynomial.  The file is comma-delimited, with the first
slot containing the HOMFLY-PT output from \texttt{jhomfly}.  This is
followed by a comma-separated list of knot types (prime and composite
through crossing number 16) matching the given HOMFLY-PT polynomial.

The prime knot types have a prefix of ``\texttt{p}'', ``\texttt{m}'',
``\texttt{h}'', or ``\texttt{a}''.  The prefix ``\texttt{a}'' means that the knot
type is achiral.  The prefix ``\texttt{h}'' means that the knot type is
chiral, but both chiralities (i.e., a configuration realizing the knot type and its
mirror image) share the same HOMFLY-PT polynomial.  If
the knot type is chiral and the HOMFLY-PT polynomial differs between
the two chiralities, then one of the chiralities has been assigned ``\texttt{p}''
(denoting ``plus'') and the other had been assigned ``\texttt{m}'' (denoting
``minus'') using the process described in the next paragraph.  Note
that the ``\texttt{p}'' versus ``\texttt{m}'' designations were assigned by
author Rawdon.  They are not ``the standard'', yet anybody performing
such a designation is likely to agree with these for most of the knot
types.  If one disagrees with the designations, just know that the 
``\texttt{p}'' versus ``\texttt{m}'' versions of a knot type are mirror images.

To designate a ``\texttt{p}'' versus ``\texttt{m}'' for a chiral knot type
whose mirror images have different HOMFLY-PT polynomials, we started by
computing the planar writhe of a minimum-crossing diagram
representative of the knot type.  If the value was positive, then the
knot type was assigned a ``\texttt{p}'', and its mirror image ``\texttt{m}''.
If the value was negative, then the knot type was assigned an
``\texttt{m}'', and its mirror image ``\texttt{p}''.  If the writhe of the
diagram was zero, then we performed knot energy minimization in
KnotPlot on a configuration forming the given knot type.  The
assignment of ``\texttt{p}'' versus ``\texttt{m}'' was then based on the sign
of the spatial writhe of the minimizing configuration.
Note that some of these spatial writhes were very close to zero, and
this is where our designations could disagree from other strategies.

For prime knot types through crossing number 10, we use ABC notation
(described in Section \ref{sec:peculiarities}).  For prime knot types
from crossing number 11 through 16, we use a slight variation of DT
notation (again, see Section \ref{sec:peculiarities}).  For example,
the knot type \texttt{-11a99} is denoted \texttt{m11.99a} and the knot
type \texttt{+11n27} is denoted \texttt{p11.27n}.

For composite knot types, the prime factor knots are ordered by crossing
number and factor knot types are separated by \#, which denotes the
connected sum (i.e., knot composition).  If there are two chiralities
of the same base knot type appearing as factors, then the \texttt{p}
factor(s) appear before the \texttt{m} factor(s).  For example,
$\texttt{p3.1\#m3.1}$ is in the table, but not $\texttt{m3.1\#p3.1}$.

\section*{Acknowledgements}
  We would like to thank Xin Liu, Renzo Ricca, and Beijing University
  of Technology for hosting the conference on Knotted Fields and for
  their enthusiasm in this endeavor.
  RGS would like to thank Chris Soteros at the University of
  Saskatchewan for support in developing recent versions of
  \texttt{jhomfly} and \texttt{jidknot} as well as Mariel V\'azquez at
  the University of California Davis for support during the original
  development of those tools, as well as \texttt{xinger}.  This
  material is based upon work supported by the National Science
  Foundation under Grant No.~1720342 to EJR.  EJR also wishes to thank
  Addie McCurdy, whose work was critical in creating the tables which
  were used to count chiral knot types.  EJR also thanks the
  University of St.~Thomas's Center for Applied Mathematics, which has
  supported several student projects related to this work.



\bibliographystyle{alpha}
\bibliography{knots}

\end{document}